\documentclass[12pt]{amsart}
\usepackage{amsthm}


\begin{document}

\def\contentsname{Table of Contents}
\def\listfigurename{List of Figures}
\def\listtablename{List of Tables}
\def\bibname{Bibliography}
\def\refname{References}
\def\indexname{Index}
\def\figurename{Figure}
\def\tablename{Table}
\def\chaptername{Chapter}
\def\appendixname{Appendice}
\def\partname{Part}
\def\abstractname{Abstract}
\newtheorem{aaa}{}[section]
\newtheorem{defi}[aaa]{Definition}
\newtheorem{lem}[aaa]{Lemma}
\newtheorem{prop}[aaa]{Proposition}
\newtheorem{cor}[aaa]{Corollary}
\newtheorem{thm}[aaa]{Theorem}
\hoffset=-1cm
\newcommand{\archi}[2]{{\footnotesize{#1~(#2)}}}
\def\noi{\noindent }
\def\noib{\noindent $\bullet $~}	
\def\pointir{\discretionary{.}{}{.\kern.35em---\kern.7em}}
\def\ie{{\textsl i.e.~}}
\def\eg{{\textsl e.g.~}}
\def\resp{{\textsl resp.~}}
\def\C{\mbox{l\hspace{-0.47em}C}}
\def\R{\mbox{I\hspace{-0.15em}R}}
\def\N{\mbox{I\hspace{-0.15em}N}}
\def\Z{\mbox{Z\hspace{-0.3em}Z}}
\def\K{\mbox{I\hspace{-0.15em}K}}
\def\E{\mbox{I\hspace{-0.15em}E}}
\def\H{\mbox{I\hspace{-0.15em}H}}
\def\Q{\mbox{l\hspace{-0.47em}Q}}
\newcommand{\cA}{{\mathcal A}}
\newcommand{\cB}{{\mathcal B}}
\newcommand{\cC}{{\mathcal C}}
\newcommand{\cD}{{\mathcal D}}
\newcommand{\cE}{{\mathcal E}}
\newcommand{\cF}{{\mathcal F}}
\newcommand{\cG}{{\mathcal G}}
\newcommand{\cH}{{\mathcal H}}
\newcommand{\cI}{{\mathcal I}}
\newcommand{\cJ}{{\mathcal J}}
\newcommand{\cK}{{\mathcal K}}
\newcommand{\cL}{{\mathcal L}}
\newcommand{\cM}{{\mathcal M}}
\newcommand{\cN}{{\mathcal N}}
\newcommand{\cP}{{\mathcal P}}
\newcommand{\cR}{{\mathcal R}}
\newcommand{\cS}{{\mathcal S}}
\newcommand{\cT}{{\mathcal T}}
\newcommand{\cV}{{\mathcal V}}
\newcommand{\cX}{{\mathcal X}}
\newcommand{\cY}{{\mathcal Y}}
\newcommand{\cZ}{{\mathcal Z}}
\def\bra{\langle }
\def\ket{\rangle }
\newcommand{\ip}[2]{\langle #1 , #2\rangle}
\def\aco{\lbrace }
\def\acf{\rbrace }
\def\pf{\par{\textsl \noindent Proof}\pointir}
\def\qed{\quad\hbox{\hskip 1pt\vrule width 4pt height 4pt
            depth 1.5pt\hskip 1pt}}
\def\qedd{\hfill\hbox{\hskip 1pt\vrule width 4pt height 4pt
            depth 1.5pt\hskip 1pt}}            
 
\font\bbbld=msbm10 scaled\magstep1
\newcommand{\bfR}{\hbox{\bbbld R}}
\newcommand{\bfC}{\hbox{\bbbld C}}
\newcommand{\bfZ}{\hbox{\bbbld Z}}
\newcommand{\bfH}{\hbox{\bbbld H}}
\newcommand{\bfQ}{\hbox{\bbbld Q}}
\newcommand{\bfN}{\hbox{\bbbld N}}
\newcommand{\bfP}{\hbox{\bbbld P}}
\newcommand{\bfT}{\hbox{\bbbld T}}
\def\Sym{\mathop{\rm Sym}}
\newcommand{\halo}[1]{\Int(#1)}
\def\Int{\mathop{\rm Int}}
\def\Re{\mathop{\rm Re}}
\def\Im{\mathop{\rm Im}}
\newcommand{\union}{\cup}
\newcommand{\goesto}{\rightarrow}
\newcommand{\bdy}{\partial}
\newcommand{\n}{\noindent}
\newcommand{\p}{\hspace*{\parindent}}

\newtheorem{theorem}{Theorem}[section]
\newtheorem{assertion}{Assertion}[section]
\newtheorem{proposition}{Proposition}[section]
\newtheorem{lemma}{Lemma}[section]
\newtheorem{definition}{Definition}[section]
\newtheorem{claim}{Claim}[section]
\newtheorem{corollary}{Corollary}[section]
\newtheorem{observation}{Observation}[section]
\newtheorem{conjecture}{Conjecture}[section]
\newtheorem{question}{Question}[section]
\newtheorem{example}{Example}[section]

\newenvironment{remark}{\smallskip\noindent{\bf Remark.}\hskip \labelsep}%

\setlength{\baselineskip}{1.2\baselineskip}

\begin{abstract}
Motivated by recent interest in the spectrum of the Laplacian of 
incomplete surfaces with isolated conical singularities, we 
consider more general incomplete 
$m$-dimensional manifolds with singularities on sets of codimension at 
least $2$.   With certain restrictions on the metric, we establish that 
the spectrum is discrete and satisfies Weyl's asymptotic formula.  
\end{abstract}

\date{}

\title[ Spectra for incomplete manifolds]{Discrete spectrum and 
Weyl's asymptotic formula for incomplete manifolds}
\author{Jun Masamune and Wayne Rossman}

\maketitle

\section{Discreteness of the Spectrum}

When one studies the Morse index of minimal surfaces in Euclidean 3-space 
$\bfR^3$ or of mean curvature $1$ surfaces in hyperbolic 3-space $\bfH^3$, the 
problem reduces to the study of the number of eigenvalues less than $2$ of 
the spectrum of the Laplace-Beltrami 
operator on Met$_1$ surfaces \cite{FC}, \cite{UY}, \cite{LR}.  
(Met$_1$ surfaces are incomplete $2$-dimensional manifolds 
with constant curvature $1$ and isolated conical singularities.)  
Met$_1$ surfaces are known to have 
pure point spectrum and satisfy Weyl's asymptotic formula.  

Here we will show that the spectrum is discrete and that Weyl's asymptotic 
formula holds for more general incomplete manifolds.  We 
allow the dimension to be arbitrary; we do not make any specific assumptions 
about the curvature; and we allow more general singularities, of at least 
codimension $2$ (in a sense to be made precise below).  This more general 
setting allows us to consider singularities such as a product of an 
$m-n$ dimensional metric cone with a portion of 
$\bfR^n$ ($m \geq n+2$), one of our 
desired examples.  In this example, the incomplete metric is 
singular only in the direction of the metric cone and not on the portion of 
$\bfR^n$ itself, so generally the incomplete 
manifolds and their metrics $\tilde{g}$ that we consider 
will not be conformally 
equivalent to open sets of compact Riemann manifolds, unlike the case 
of Met$_1$ surfaces.  With this in mind, we now define the types of incomplete 
manifolds and metrics $\tilde{g}$ that we will study here.  

Let $(M,g)$ be a compact manifold of dimension $m$ with 
smooth Riemannian metric 
$g$.  Let $N$ be a compact submanifold of dimension $n$ with codimension 
$m-n \geq 2$.  Suppose further that in a neighborhood of $N$ the metric 
$g$ can be diagonalized; that is, there exist 
local coordinates $(x_1,...,x_{m-n},y_1,...,y_n)$, where 
$(0,...,0,y_1,...,y_n)$ are coordinates for $N$, so that 
$(dx_1,...,dx_{m-n},dy_1,...,dy_n)$ is globally defined in some open 
neighborhood of $N$ and so that the metric $g$ is written 
\[ g = \left( \begin{array}{cccc}
g_1 & 0 \\
0 & g_2 
\end{array} \right) \; , \] where $g_1$ is an $m-n \times m-n$ positive 
definite matrix, and $g_2$ is an $n \times n$ positive definite matrix.  (For 
example, such a 
case can occur if $M$ has a product structure $M = M_1 \times N$ near 
$N$, where 
$M_1$ is an $m-n$ dimensional compact Riemannian manifold.)  

\begin{theorem} \label{firstresult} 
Let $N$ be an $n$-dimensional compact submanifold of an $m$-dimensional 
compact manifold $(M,g)$ with $m \geq n+2$ such that the metric $g$ can be 
diagonalized near $N$.  Choose local coordinates in a neighborhood of 
$N$ so that 
\[ g = \left( \begin{array}{cccc}
g_1 & 0 \\
0 & g_2 
\end{array} \right) \] in this neighborhood.  
Let $\tilde{g}$ be another smooth regular metric on $M \setminus N$ so that 
\[ \tilde{g} = \left( \begin{array}{cccc}
f^2 g_1 & 0 \\
0 & g_2 
\end{array} \right) \] in a neighborhood of $N$, where 
$f \in C^\infty(M \setminus N)$.  

If $m=2$, assume that $f \in L_g^{2+\epsilon}(M)$ for some $\epsilon \in 
(0,\infty)$.  

If $m \geq 3$, assume that $\inf(f)>0$ and 
$f \in L_g^{(m(m-n)/2) +\epsilon}(M)$ for some 
$\epsilon \in (0,\infty)$.  

Then 
the Sobolev space $W^{1,2}_{\tilde{g}}(M \setminus N)$ with respect to 
$\tilde{g}$ is compactly included in 
$L^2_{\tilde{g}}(M \setminus N)$.  
\end{theorem}

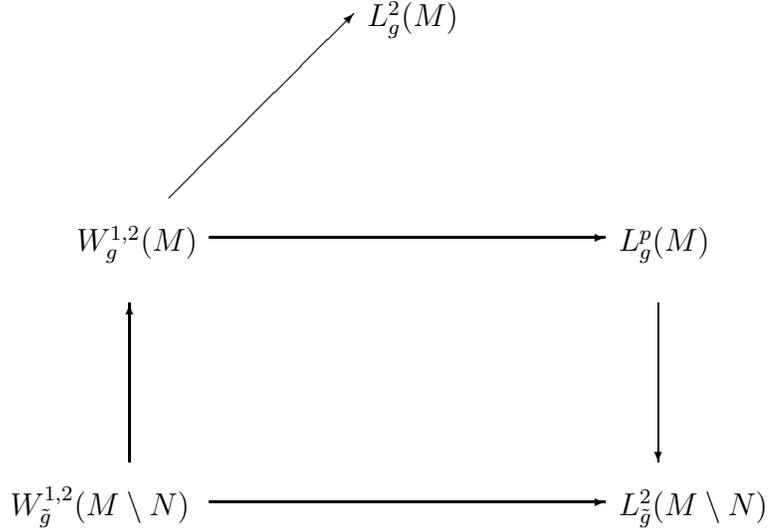
\begin{figure}
\begin{picture}(300,200)
\put(20,0){$W^{1,2}_{\tilde{g}}(M \setminus N)$}
\put(45,100){$W^{1,2}_g(M)$}
\put(250,0){$L^2_{\tilde{g}}(M \setminus N)$}
\put(250,100){$L^p_g(M)$}
\put(155,185){$L^2_g(M)$}
\put(65,20){\vector(0,1){60}}
\put(95,5){\vector(1,0){150}}
\put(95,105){\vector(1,0){150}}
\put(265,80){\vector(0,-1){60}}
\put(80,120){\vector(1,1){70}}
\end{picture}
\caption{The compact inclusion of $W^{1,2}_{\tilde{g}}(M \setminus N)$ into 
$L^2_{\tilde{g}}(M \setminus N)$.}  
\end{figure}

\begin{proof}
When $m \geq 3$ and $p \in (2,2m/(m-2))$ (resp. $m=2$ and $p \in (2,\infty)$), 
then the inclusion $W^{1,2}_g(M)$ into $L_g^p(M)$ is compact.  
When $m \geq 3$ and $f \in L^{(m(m-n)/2) + \epsilon}$ (resp. 
$m =2$ and $f \in L^{2 + \epsilon}$) for some 
positive $\epsilon$, then the inclusion $L^p_g(M)$ into 
$L^2_{\tilde{g}}(M \setminus N)$ is continuous, 
by H\"older's inequality.  
For example, when $m \geq 3$, we can choose 
\[ p = {m + (2 \epsilon / (m-n)) 
\over (m / 2) + (\epsilon / (m-n)) - 1} \; , \] and then the 
H\"older inequality implies 
\[ \|u\|_{L^2_{\tilde{g}}} = \sqrt{\int u^2 f^{m-n} dA} \leq c \cdot 
\|u\|_{L_g^p} \] for \[ 
c = \left(\int f^{(m (m-n)/2)+\epsilon} dA 
\right)^{((m / 2) + (\epsilon / (m-n)))^{-1}/2} < \infty \; . \] 
So we only need to show that $W^{1,2}_{\tilde{g}}(M \setminus N)$
is continuously contained in $W^{1,2}_g(M)$ to conclude
$W^{1,2}_{\tilde{g}}(M  \setminus N)$ is compactly contained in
$L^2_{\tilde{g}}(M  \setminus N)$.  When $m \geq 3$, this is clear, 
since $\inf(f) > 0$.  
When $m=2$, then $n=0$, and $g$ and $\tilde{g}$ are conformally
equivalent on $M \setminus N$.  
Suppose by way of contradiction that the inclusion is not continuous, that is, 
that there exists a sequence of functions $u_k$ such 
that $\|u_k\|_{W^{1,2}_{g}} = 1$ and $\|u_k\|_{W^{1,2}_{\tilde{g}}}
 < 1 / k$.  By choosing a subsequence if necessary,
we may assume the following:
\begin{enumerate}
\item there exists a function $u$ such that $u_k \to u$, $W^{1,2}_g$-weakly,
\item there exists a function $v$ such that $u_k \to v$, $L^p_g$-strongly,
\item $u_k \to v$, $L^2_{\tilde{g}}$-strongly,
\item $u_k \to v$, $L^2_g$-strongly.  
\end{enumerate}
The fourth item follows from the fact that $\|u_k - v\|_{L^2_g} \leq 
\hat{c} \cdot \|u_k - v\|_{L^p_g}$, since $(M,g)$ is smooth and compact.  
As $u_k$ converges to both $u$ and $v$ $L^2_g$-weakly, $u=v$.  Also, 
\[ 1 = \liminf_{k \to \infty} \|u_k\|_{W^{1,2}_g} \geq
 \|u\|_{W^{1,2}_g} \; . \]
Let $\nabla$ and $dA$ (resp. $\tilde{\nabla}$ and $d\tilde{A}$) denote the 
gradient and area-form with respect to the metric $g$ (resp. $\tilde{g}$).  
Then, using $\int_M |\nabla u_k |_g^2 dA = \int_M | \tilde{\nabla} u_k 
|_{\tilde{g}}^2 d\tilde{A}$, we have $\int_M u_k^2 dA \to 1$ and 
$\int_M u^2 dA = 1$ and $\int_M | \nabla u |_g^2 dA = 0$, so 
$u$ is a nonzero constant.  Also, $\int_M u_k^2d\tilde{A} \to \int_M u^2 
d\tilde{A} = 0$, so $\int_M d\tilde{A} = 0$.  This is a contradiction, since 
$f$ is smooth on $M \setminus N$ and not identically zero.  
\end{proof}

\begin{remark}
For $m \geq 3$, the condition $\inf(f) > 0$ is a simple way to ensure 
$W_{\tilde{g}}^{1,2}$ is continuously contained in $W_g^{1,2}$, but it is 
necessary.  This is not generally a continuous inclusion if $\inf(f)=0$.  
For example, suppose $\inf(f) = 0$, and $n=0$.  Let 
$M_k = \{p \in M \setminus N \,\, | \,\, |f(p)| < 1 / k \}
 \neq \emptyset$.  
Choose $u_k$ so that supp$(u_k) \subset M_k$ and 
$\|u_k\|_{W_g^{1,2}}^2 = 1$.  Then $\tilde{g} = f^2 g$ and $g$ are
conformally equivalent and 
\[ \|u_k\|_{W_{\tilde{g}}^{1,2}}^2 = \int_{M_k} (u_k)^2 f^m dA + 
\int_{M_k} | \nabla u_k |_g^2 f^{m-2} dA \leq {1 \over k^{m-2}} 
\|u_k\|_{W_g^{1,2}}^2 = {1 \over k^{m-2}} \to 0 \] as 
$k \to \infty$.  Hence, we do not have continuous inclusion.  
\end{remark}

\medskip

Let $\overline{\triangle}_{\tilde{g}}^{F}$ denote the Freidrichs' self-adjoint 
extension of the Laplacian with domain $C_0^\infty(M \setminus N)$, and let 
$W^{1,2}_{0,\tilde{g}}(M \setminus N)$ be the closure of 
$C_0^\infty(M \setminus N)$ in the $W^{1,2}_{\tilde{g}}(M \setminus 
N)$ norm.  
Standard arguments give the following: 

\begin{corollary}
Let $(M \setminus N,\tilde{g})$ be as in Theorem {\em \ref{firstresult}}.  
The operator 
$\overline{\triangle}_{\tilde{g}}^{F}$ on 
$(M \setminus N,\tilde{g})$ has discrete spectrum consisting of eigenvalues 
$0 = \lambda_1 < \lambda_2 \leq \dots \leq \lambda_j \leq \dots \to 
+ \infty$, each with 
multiplicity $1$.  The corresponding eigenfunctions $\phi_1,\phi_2,... \in 
W^{1,2}_{\tilde{g}}(M \setminus N)$ can be chosen as an orthonormal 
basis for $L^{2}_{\tilde{g}}(M \setminus N)$.  
Furthermore, the variational characterization for the eigenvalues holds:  
\[ \lambda_j = \inf_{V^j} \sup_{\phi \in V^j, \phi \neq 0} 
\frac{\| \nabla \phi \|_{L^{2}_{\tilde{g}}(M \setminus N)}^2}
{\| \phi \|_{L^{2}_{\tilde{g}}(M \setminus N)}^2} \; , \] where 
$V^j$ represents an arbitrary $j$-dimensional subspace of 
$W^{1,2}_{0,\tilde{g}}(M \setminus N)$.  
\end{corollary}

\begin{remark}
When 
$W^{1,2}_{0,\tilde{g}}(M \setminus N) = W^{1,2}_{\tilde{g}}(M \setminus N)$, 
$(M \setminus N,\tilde{g})$ has negligible boundary,
in Gaffney's sense \cite{G}.   
Therefore, the Laplacian considered on Gaffney's domain of functions is 
essentially self-adjoint.  Let $\overline{\triangle}_{\tilde{g}}^{G}$
denote the unique self-adjoint extension.  One can show 
that the two self-adjoint operators 
$\overline{\triangle}_{\tilde{g}}^{F}$ and 
$\overline{\triangle}_{\tilde{g}}^{G}$ have equal domains, that is, 
\[ \triangle := \overline{\triangle}_{\tilde{g}}^{G} = 
\overline{\triangle}_{\tilde{g}}^{F}
\; . \]  Since $(M \setminus N,\tilde{g})$ has negligible boundary, this operator's 
domain has no boundary conditions at $N$.  
So when 
$W^{1,2}_{0,\tilde{g}}(M \setminus N) = W^{1,2}_{\tilde{g}}(M \setminus N)$, 
this is the operator for which we will study the spectrum, and it is the 
same operator as that used in the study of Morse index of minimal surfaces in 
$\bfR^3$ and mean curvature 1 surfaces in $\bfH^3$.  
\end{remark}

\medskip

As seen in the above remark, we would like to consider the cases where 
$W^{1,2}_{0,\tilde{g}}(M \setminus N) = W^{1,2}_{\tilde{g}}(M \setminus N)$.  
We will also need this property for establishing Weyl's asymptotic formula, 
so we now give a sufficient condition to imply this property \cite{M2}.
In order to state it, here we introduce the notion of capacity and
Cauchy boundary \cite{M1},\cite{M2}.

\begin{definition}
   Let $M$ be an arbitrary Riemannian manifold.
   We denote by $\mathcal{O}$, the family
   of all open subsets of the completion $\overline{M}$ of $M$.
   For $A\in\mathcal{O}$, we define the set of functions $L_{A}$ by
   \[
    L_{A}
    =\{ f\in W^{1,2}(M) \,|\, f \ge 1 \ \mbox{a.e. on}\ A \}.\]
    We define the capacity of $A$, $\mathrm{Cap}(A)$, by
   \begin{eqnarray*}
           & & \mathrm{Cap}(A)=\left\{
               \begin{array}{ll}
               \displaystyle
               \inf_{f \in L_{A}} \| f \|_{W^{1,2}},
                 & L_{A} \neq\phi, \\
               \infty, & L_{A}=\phi. \\
               \end{array}\right. \\
   \end{eqnarray*}
   For a Borel set $B\subset\overline{M}$,
   we define the capacity $\mathrm{Cap}(B)$ by
   \[\mathrm{Cap}(B)=\displaystyle\inf_{A\in\mathcal{O},\,
               \mathit{B}\subset \mathit{A}}
               \mathrm{Cap}(A).\]
   We say that a subset $B$ of $\overline{M}$ is 
   \textit{almost polar} if $\mathrm{Cap}($B$)=0$.
\end{definition}

\begin{definition}
The Cauchy boundary $\partial M$ of $M$ is defined by 
\[ \partial M := \overline{M}\setminus M,\]
where $\overline{M}$ is the completion of $M$ with respect to the Riemannian
distance.
\end{definition}

\begin{lemma}
{\em (\cite{M2})} 
For an arbitrary Riemannian manifold $M$, let $\partial M$ denote
the Cauchy boundary of $M$.
Then the two Sobolev spaces $W_{0}^{1,2}(M)$ and $W^{1,2}(M)$ coincide
if and only if $\partial M$ is an almost polar set.
\end{lemma}
In the case of Theorem 1.1, the Cauchy boundary of $M\setminus N$
is $N$.
It is shown in \cite{M1} that when the lower Minkowski codimension of the 
Cauchy boundary is not less than 2,
then $\partial M$ is almost polar,
where the lower Minkowski codimension is defined as follows:
\begin{definition}
The lower Minkowski codimension of $\partial M$  is defined to be 
\[ \underline{\mathrm{codim}}_{\mathcal M} (\partial M) :=
 \lim_{R \to 0} \inf  {\log(\mathrm{vol}({\mathcal N}_R)) \over 
\log(R)} \; , \] 
where ${\mathcal N}_R$ is a radius $R$ tubular neighborhood of $\partial M$.  
\end{definition}

We now consider some examples.  

\begin{example} {\em 
Consider the ``football''.  Set $M = \bfC \cup 
\{ \infty \}$ and $N = \{0,\infty\}$ ($m=2$ and $n=0$) and set 
\[ g = {4(dx^2+dy^2) \over(1+r^2)^2} \; , \; \; \; 
f = {\mu r^{\mu-1} (1+r^2) \over 1+r^{2\mu}} \; , \; \; \; \mu \in \bfR^+ 
\; , \; \; \; \tilde{g} = f^2 g \; , \] 
where $r = \sqrt{x^2+y^2}$.  Note that $f \in L_g^{2+\epsilon}(M \setminus 
N)$ for some $\epsilon > 0$, and 
$\underline{\mbox{codim}}_{\mathcal M} (N) = 2$ for any $\mu$.  
When $\mu < 1$, the football is an Alexandrov 
space and $\triangle$ has discrete spectrum, by \cite{KMS} or by Theorem 
\ref{firstresult} above.  
When $\mu > 1$, the football is not an Alexandrov space, but the 
spectrum is still discrete, by Theorem \ref{firstresult} (see also Lemma 
4.3 of \cite{LR}).  }
\end{example}

\begin{example} {\em 
Consider a compact $m$-dimensional manifold $M$ with 
metric $g$, and suppose $M$ contains the unit ball $B^m$ so that 
$g$ is the standard Euclidean metric on $B^m \subset M$.  
Let $N= \{ \vec{o} \}$ be the center point of $B^m \subset M$.  Let 
$f=r^\ell$ on $B^m$ and $\tilde{g} = f^2 g$ with 
$\ell \in (-2 / m,0)$, and extend $f$ to be positive and smooth 
on $M \setminus B^m$.  
Thus $(M \setminus N,\tilde{g})$ is not complete, and 
$\underline{\mbox{codim}}_{\mathcal M} (N) = m$ for any $\ell$.  Also, as
$f$ satisfies the conditions of Theorem \ref{firstresult}, $\triangle$ on 
$(M \setminus N,\tilde{g})$ has discrete spectrum.  }
\end{example}

\begin{example} {\em 
As it is known that Alexandrov spaces have discrete spectrum 
\cite{KMS}, we are interested in finding examples that are not Alexandrov 
spaces and for which Theorem \ref{firstresult} can be applied.  The footballs 
with $\mu > 1$ provide such examples in two dimensions.  The 
following example shows that one can easily find such examples in higher 
dimensions as well.  (We choose a slightly complicated function $f$ in 
order to easily verify it will not be an Alexandrov space.)  

Consider the previous example with $m=3$; that is, $M$ is compact, 
$3$-dimensional, $N = \{ \vec{o} \} \subset B^3 \subset M$, and 
$g$ is the Euclidean metric on $B^3$.  Set $f = \cos^2(\phi) r^\ell 
+ \sin^2(\phi) (1+r^{3 \over 2}) \in L_g^{(9/2)+\epsilon}(M)$ on 
$B^3 \subset M$, where $(r,\theta,\phi)$ are the 
spherical coordinates of $B^3$, and $\ell \in (-2 / 3,0)$.  Extend 
$f$ to be positive and smooth on $M \setminus B^3$, and let 
$\tilde{g} = f^2 g$.  Then 
$(M \setminus N,\tilde{g})$ is not complete, and the conditions of Theorem 
\ref{firstresult} are satisfied.  Hence the spectrum  of 
$\overline{\triangle}_{\tilde{g}}^{F}$ is discrete.  The ball 
$B^3$ is invariant under the isometry 
$(r,\theta,\phi) \to (r,\theta,\pi-\phi)$, thus the 
sectional curvature in the $\{\phi = \pi / 2 \}$-plane is 
$K_{\tilde{g}} = -(\triangle \ln(f) ) 
/ f^2 \to - \infty$ near $N$, so it is not an Alexandrov space.  }
\end{example}

\begin{example} {\em 
Consider the $3$-dimensional torus $M = T^3 = \bfR^3 / \bfZ^3$ with the 
standard Euclidean metric $g$, and 
the 1-dimensional torus $N = S^1 = (\bfR / \bfZ, 0, 0) \subset M$.  
We will use cylindrical coordinates $(x,r,\theta)$, where 
$r$ is the radial distance to $N$ and $x$ is the arc-length along $N$.  
Let $f = \cos^2(\theta) + \sin^2(\theta) r^\ell$ near $N$ with 
$\ell \in (-2 / 3,0)$, and extend $f$ to be positive and smooth 
away from $N$.  This manifold is incomplete, and 
$\overline{\triangle}_{\tilde{g}}^{F}$ has discrete spectrum, by Theorem 
\ref{firstresult} and Corollary 1.1.  }
\end{example}

\begin{remark} 
Suppose $M$ is $2$-dimensional and contains $B^2$ so that $g$ is the standard 
Euclidean metric when restricted to $B^2$.  Suppose 
$N = \{\vec{o}\} \subset B^2 \subset M$ and 
$f = -1 / (r \ln(r))$ near $N$.  Then, with respect to $\tilde{g} = 
f^2 g$, we have a complete end at $N$ that is a  
curvature $-1$ psuedosphere of finite area, so 
the spectrum is not discrete \cite{D}, \cite{Mu}.  
Since $f \in L^2(M,g)$, but $f \not\in L^{2+\epsilon}(M,g)$ for all positive 
$\epsilon$, we know Theorem \ref{firstresult} is sharp when $m=2$.  
(If we had chosen $f = 1 / r \in L^{2-\epsilon}(M,g)$ for all small 
positive $\epsilon$ instead, we would have produced a 
round cylindrical end of radius $1$ which does not have discrete spectrum and 
does not have finite area.)  
\end{remark}

\begin{remark}
Consider $M = T^2 \times T^{m-2}$ and $N = T^{m-2}$ and 
$f = -1 / (r \ln(r))$ near $N$, where $r$ is radial distance to $N$.  
Let the diagonalized coordinates near $N$ be $(x_1,x_2,y_1,...,y_{m-2})$, 
inherited from the standard rectangular Euclidean coordinates of $\bfR^m$.  
Then $(M \setminus N,\tilde{g})$ is complete, and the sectional curvatures are 
\[ K_{\tilde{g}}(\partial_{x_1},\partial_{x_2}) = -1 \; , \; \; 
K_{\tilde{g}}(\partial_{x_i},\partial_{y_j}) = 0 \; , \; \; 
K_{\tilde{g}}(\partial_{y_i},\partial_{y_j}) = 0 \; . \] 
So the Ricci curvature is bounded 
below, and hence the essential spectrum is not empty \cite[Theorem 3.1]{D}.  
So Theorem \ref{firstresult} is not 
true for this $f \in L^p$, $p \leq 2$.  Hence, for all 
$m$, the restriction on $f$ in 
Theorem \ref{firstresult} cannot be weakened to $f \in L^p$
for some $p \leq 2$.  
\end{remark}

\begin{remark}
Donnelly and Li \cite{DL} have found complete examples
 $(M \setminus N,\tilde{g})$ where 
$M = \bfR^m \cup \{ \infty \}$ and $N = \infty$ ($n=0$) and $\tilde{g}$ is 
rotationally invariant, so that sectional curvature converges to 
$-\infty$ at $N$ and $(M \setminus N,\tilde{g})$ has 
pure point spectrum.  For example, let $m=2$ and 
$\tilde{g} = dr^2+e^{-r^k} d\theta^2, \, k>1$, in radial coordinates
 $(r,\theta)$ of 
$\bfR^2$.  It is complete 
and its single end is conformally a punctured disk, and since
 the curvature converges to 
$-\infty$ at the end, it has pure point spectrum \cite{DL}.  Theorem 
\ref{firstresult} does not apply to such examples.  
\end{remark}

\section{Weyl's formula}

In this section, let $M$ be an $m$-dimensional Riemannian manifold
with finite volume and finite diameter. $M$ can be noncompact and 
incomplete.  

\begin{remark}
The manifolds $(M \setminus N,\tilde{g})$ in Theorem 1.1 
have finite volume, since 
$f \in L_g^{(m (m-n)/2) + \epsilon} \subseteq L_g^{m-n}$ implies 
$\mbox{vol}(M \setminus N,\tilde{g}) = \int_M d\tilde{A} = 
\int_M f^{m-n} dA < \infty$.  
\end{remark}

\medskip

Before stating and proving Weyl's asymptotic formula, we 
establish some notation.  
Let ${\mathcal N}_R$ be a 
radius $R$ tubular neighborhood of the Cauchy boundary
$\partial M$ of $M$.  Note that 
vol$(M \setminus {\mathcal N}_R)+$vol$({\mathcal N}_R)=$vol$(M)$ and 
vol$({\mathcal N}_R) \to 0$ as $R \to 0$.  
Define the Neumann isoparimetric constant of ${\mathcal N}_R$ by 
\[ C_R := \inf_\gamma \frac{\mbox{vol}(\gamma)}
{\min \{ \mbox{vol}(M_1), \mbox{vol}(M_2) \}^{(m-1) / m} } \; , \] 
where the infimum is taken over all hypersurfaces $\gamma$ of 
${\mathcal N}_R$ which divide ${\mathcal N}_R$ into two parts $M_1$ and 
$M_2$, and 
where $\mbox{vol}(\gamma)$ represents the $m-1$ dimensional volume of 
$\gamma$ and $\mbox{vol}(M_j)$ represents the $m$-dimensional
volume of $M_j$.  

Here, we will assume that 
\[ C := \inf_{R > 0} C_R > 0 \; . \]
Then, since $M$ has finite 
volume, one can see that $M$, $N_{R}$ and $M\setminus N_{R}$ all 
have pure point spectra.  
Let $\lambda_j^{1,N}$ (resp. $\lambda_j^{2,N}$) be the Neumann eigenvalues 
on $\mbox{Int}(M \setminus {\mathcal N}_R)$ (resp. 
$\mbox{Int}({\mathcal N}_R) $) counted with their multiplicities (i.e. 
listed in nondecreasing order, and the number of times that 
any eigenvalue appears in the list equals its multiplicity).  Let 
$\lambda_j^{3,N}$ be the Neumann eigenvalues of 
$\mbox{Int}(M \setminus {\mathcal N}_R) \cup \mbox{Int}({\mathcal N}_R)$ 
counted with their multiplicities.  Let $\lambda_j^{D}$ be the Dirichlet 
eigenvalues on $\mbox{Int}(M \setminus {\mathcal N}_R)$ counted with their 
multiplicities.  Here, we state Weyl's asymptotic formula for $M$.  

\begin{theorem}
Let $M$ be an $m$-dimensional Riemannian manifold with
finite volume and finite diameter.
If the Cauchy boundary of $M$ is an almost polar set and
$C > 0$, then the eigenvalues $\lambda_j$ of the Laplacian $\triangle$
satisfy Weyl's asymptotic 
formula \[ \lim_{j \to \infty} \frac{\lambda_j^{m / 2} 
\mathrm{vol}(M)}{j} = \frac{(2 \pi)^m}{\mathrm{vol}(B^m)} \; . \] 
\end{theorem}

\begin{proof}
Let $W = (2 \pi)^m / \mathrm{vol}(B^m)$.  
Note that $\lambda_j \leq \lambda_j^D$ by Dirichlet-Neumann 
bracketing techniques (see, for example, volume 4 of \cite{RS}).  Note 
also that, since $M$ has finite diameter and therefore 
$M \setminus {\mathcal N}_R$ is relatively compact, 
the $\lambda_j^D$ satisfy Weyl's asymptotic formula on
$M \setminus {\mathcal N}_R$.  So 
\[ \lambda_j \leq \lambda_j^D \approx W \mbox{vol}(M \setminus 
{\mathcal N}_R)^{-2 / m} j^{2 / m} \to 
W \mbox{vol}(M)^{-2 / m} j^{2 / m} \] for large $j$, 
as $R \to 0$.  This implies 
\[ \limsup_{j \to \infty} \frac{\lambda_j^{m / 2} 
\mbox{vol}(M)}{j} \leq W \; . \] 

Consider the Neumann heat kernel 
\[ H_R(x,y,t) = \sum_{i=1}^\infty e^{-\lambda_i^{2,N} t} \phi_i^{2,N}(x) 
\phi_i^{2,N}(y) \] on ${\mathcal N}_R$, where $\{ \phi_i^{2,N} \}_{i=1}^\infty$
is an orthonormal basis of eigenfunctions in $L^2({\mathcal N}_R)$
associated to the eigenvalues $\lambda_i^{2,N}$.  
Using the method of \cite{LT}, we know that the Neumann heat kernel on 
${\mathcal N}_R$ belongs to the Sobolev space $W^{1,2}({\mathcal N}_R)$ and 
has the above form.  
As the isoperimetric constant $C_R$ of ${\mathcal N}_R$ is positive and the 
coarea formula on ${\mathcal N}_R$ holds for nonnegative functions, the 
associated Neumann Sobolev constant of 
${\mathcal N}_R$ is also positive.  Additionally, we have $H_R(x,y,t)$ in 
the above form, so the methods in \cite{CL} can be applied to show  
\[ \lambda_j^{2,N} \geq \alpha(m) C_R^2 \left({j \over 
\mbox{vol}({\mathcal N}_R)}\right)^{2 / m} \geq \alpha(m)
C^2 \left({j \over 
\mbox{vol}({\mathcal N}_R)}\right)^{2 / m} \; , \] where 
$\alpha(m)$ is a positive constant depending only on $m$.  

Note that the list $\{ \lambda_j^{3,N} \}$ is equal to the disjoint union 
of the lists $\{ \lambda_j^{1,N} \}$ and $\{ \lambda_j^{2,N} \}$ rearranged in 
increasing order.  Note also that $\lambda_j \geq \lambda_j^{3,N}$, by 
Dirichlet-Neumann bracketing.
 Since $\lambda_j^{2,N} \geq \alpha(m) C^2 (j /
\mbox{vol}({\mathcal N}_R))^{2 / m}$ and $\lambda_j^{1,N} \approx 
W \mbox{vol}(M \setminus {\mathcal N}_R)^{-2 / m} j^{2 / m}$ for large $j$, 
and since vol$({\mathcal N}_R) \to 0$ and 
vol$(M \setminus {\mathcal N}_R) \to$vol$(M)$ as $R \to 0$, we have 
\[ \liminf_{j \to \infty} \frac{\lambda_j^{m / 2} 
\mbox{vol}(M)}{j} \geq W \; . \] 
\end{proof}

\begin{example} {\em Examples 1.1 and 1.2 satisfy the conditions of
Theorem 2.1, hence their 
eigenvalues satisfy Weyl's asymptotic formula.  }
\end{example}

\begin{remark}
Using the methods of \cite{CL}, one can additionally conclude that 
$\lambda_j^{m / 2} \geq \alpha (m) \hat{C}^{m / 2} 
j/\mbox{vol}(M)$ for some positive constant 
$\hat{C}$ depending only on the lower bound of 
the Sobolev constants of $\mathcal N_{R}$ for all $R>0$.
\end{remark}

\begin{remark}
Because the ``football'' in Example 1.1 satisfies the conditions of 
Theorem 2.1, it is clear that all Met$_{1}$ surfaces also satisfy the 
conditions of Theorem 2.1.  The authors hope to consider the more general 
case where the conical singularities
form a fractal set, and hope that Theorem 2.1 can be applied to such cases.
As an example of such a case, since Minkowski dimension
and Hausdorff dimension coincide on self-similar fractals,
the Cauchy boundary of $(S^{3}\setminus {\mathcal C}, g_{S^3})$ is
almost polar, where $\mathcal C$ is a Cantor set.
\end{remark}

\begin{remark}
The results here bear some relation to the work \cite{KS}, in which 
Kuwae and Shioya have recently studied the convergence of 
the spectra of a sequence of Riemannian manifolds (they do not 
assume completeness of the manifolds).  Some of the results in 
\cite{KS} involve the almost polarity condition.  
\end{remark}

\end{document}